\input amstex
\magnification 1200

 \documentstyle{amsppt}
\NoRunningHeads
\TagsOnRight
\pageheight {7.291in}
\pagewidth {5.5in}
\redefine \r{\langle r,r \rangle}

 \topmatter

 \title Boundary Solutions of the Classical Yang-Baxter Equation
 \endtitle
 \author M. Gerstenhaber and A. Giaquinto
 \endauthor
 \address
 Department of Mathematics, University of Pennsylvania,
 Philadelphia PA 19104-6395  \endaddress
  \email
 mgersten\@mail.sas.upenn.edu
 \endemail
 \address
 Department of Mathematics, Mississippi State University,
P.O. Drawer MA, Mississippi State MS 39762  \endaddress
 \email
 tonyg\@math.msstate.edu
 \endemail
\thanks Both authors thank the NSA for partial support of this research.
\endthanks
 \abstract
 We define a new class of unitary solutions to the classical Yang-Baxter
equation (CYBE). These ``boundary solutions'' are those which lie in the
closure of the space of unitary solutions to the modified
classical Yang-Baxter equation (MCYBE).
Using the Belavin-Drinfel'd classification
of the solutions to the MCYBE, we are able to exhibit new
families of solutions to the CYBE. In particular, using the
Cremmer-Gervais solution to the MCYBE,
we explicitly construct for all $n \geq 3$ a boundary
solution based on the maximal parabolic subalgebra of
$\frak{sl}(n)$ obtained by deleting the first negative root.
We give some evidence for a generalization of this result pertaining
to other maximal parabolic subalgebras whose omitted root is relatively
prime to $n$.
We also give
examples of non-boundary solutions for the classical simple
Lie algebras.

\endabstract
\endtopmatter
 \baselineskip 14pt

\subhead 0. Introduction
\endsubhead
Let $\frak g$ be a simple Lie algebra, either over $\Bbb C$ or
split over $\Bbb R$.  To fix the terminology call an $r \in \frak
g \bigwedge \frak g$ a (constant, ``unitary'') solution to the
{\it classical Yang-Baxter equation} (CYBE) if
$\r :=[r_{12},r_{13}] +[r_{12},r_{23}]+ [r_{13},r_{23}] = 0,$
and call it a solution to the {\it modified} CYBE (MCYBE) if
instead $\r$ is a non-zero multiple of the unique, up to
scalar multiple, invariant of $\bigwedge^3\frak g.$
The simply connected algebraic group $G$ with
Lie$\, (G)=\frak g$ operates on
both sets, which together with multiplication by non-zero scalars
defines a natural concept of equivalence.

The solutions to the MCYBE have been constructivly classified by
Belavin and Drinfel'd and depend on certain ``admissible triples''
involving the root system of $\frak g$, see [BD],
or [CP, \S 3.1.A and \S 3.1.B]. These triples have the form
$(\Pi_1, \Pi_2,T)$ where $\Pi_i$ is a proper subset of $\Pi$, the
set of positive simple roots, and $T:\Pi_1\rightarrow \Pi_2$ is
a bijection satisfying certain properties (see Definition 2.1
for a precise formulation of these properties).
These triples serve as the ``discrete''
parameter for the space of solutions to the MCYBE; there is also a
``continuous'' parameter: it is a certain
$\binom {d-\#\Pi_1}2$-dimensional affine subspace of
$\frak h \bigwedge \frak h$ where $d$ is the rank of $\frak g$
and $\frak h$ is a fixed Cartan subalgebra.
When $\frak {g}=\frak {sl}(2)$ there is only the trivial (or empty)
triple which has $\Pi_1=\Pi_2=\emptyset$ and there is a
unique solution, $e_{12}\wedge e_{21}$, to the MCYBE lying in $\frak {sl}(2)\bigwedge
\frak {sl}(2)$. For $\frak g=\frak {sl}(3)$, there are essentially
two admissible triples, the trivial triple and the
{\it Cremmer-Gervais} triple. The trivial triple determines
the following family of solutions to the MCYBE:
$$e_{12}\wedge e_{21}+e_{13}\wedge e_{31}+e_{23}\wedge e_{32} +
\lambda \cdot \left( e_{11}\wedge e_{22}-e_{11}\wedge e_{33}+
e_{22}\wedge e_{33}\right)$$
where $\lambda $ is an arbitrary scalar and is the continuous
parameter for this family. In contrast, the continuous parameter 
for the Cremmer-Gervais triple is uniquely determined. The solution
to the MCYBE associated to this triple is
$$e_{12}\wedge e_{21}+e_{13}\wedge e_{31}+e_{23}\wedge e_{32}
+\tsize {\frac {1}{3}} \left( e_{11}\wedge e_{22}-e_{11}\wedge e_{33}+
e_{22}\wedge e_{33}\right)+2\, e_{12}\wedge e_{32}.\tag 0.1$$
This Cremmer-Gervais triple (and its associated solution to the MCYBE)
have analogs in higher dimensions, (see our remarks preceeding
Theorem 2.6 and formula (5.8)).
Actually,
what we call the Cremmer-Gervais triple or Cremmer-Gervais
solution to the MCYBE is never mentioned in [CG]. Instead,
the authors there construct (a non-standard) quantum Yang-Baxter
matrix. With a knowledge of the Belavin-Drinfel'd classification,
it is easy to infer that the classical limit of their
quantum Yang-Baxter matrix is, essentially, the Cremmer-Gervais
solution mentioned above. The Cremmer-Gervais triple 
for $\frak {sl}(n)$ is described in [GGS2] (see [GGS3] for a more
detailed exposition); the case $n=3$ can
also be found in [FG1].
A derivation of the Belavin-Drinfel'd 
classification using homological and deformation-theoretic
methods can be found in [FG1] and [FG2].  

If $\frak g=\frak {sl}(n)$, there are triples such that $\Pi_1$ --  
which must omit at
least one positive simple root -- omits exactly one. This is
possible only when $n\geq 3$. Similar
considerations hold for other split simple $\frak g$ when exactly
two roots are omitted. 
Such triples determine
a unique solution to the MCYBE since the continuous parameter
has dimension zero.
Identifying the positive simple roots of $\frak {sl}(n)$ with
$\{1,2,\ldots ,n-1\}$ in their natural order,
we will prove (see Theorem 2.6) that the
index of the omitted root must be relatively prime to $n$ and $T$
is then unique:  Denoting that index by $n-i$, the root omitted
from $\Pi_2$ is then $i$ and for all $j$ between $1$ and
$n-1$ except $n-i,$ we have $T(j) = {j+i}\mod n.$
The solution of the MCYBE constructed
with such a triple will be called a {\it generalized Cremmer-Gervais} 
solution, the original Cremmer-Gervais case being that
where $i=1$ and $T(j) = j+1$ for all $1\leq j \leq n-1$. 

For solutions to the CYBE, by contrast, there is presently only a
non-constructive description due to Stolin, see [S1], [S2], or
[CP, \S 3.1.D]. A Lie subalgebra $\frak f$ of $\frak g$ is
{\it {quasi-Frobenius}} if it has a non-degenerate 2-cocycle with
coefficients in the ground field $k$; it is {\it {Frobenius}} if
it has a non-degenerate 2-coboundary, that is, if there is a
linear map $F:\frak f \to k$ such that $F([a,b])$ is a 
non-degenerate skew form on $\frak f \wedge \frak f.$  A simple Lie
algebra can not itself be Frobenius, so $\frak f$ is necessarily
a proper subalgebra. Every solution $r$ of the CYBE in $\frak g$
has a largest subalgebra $\frak g_0,$ which we will call its {\it
carrier}, on which it is non-degenerate. That subalgebra is
necessarily quasi-Frobenius, the 2-cocycle being just the inverse
of $r.$  That is, if relative to some basis $b_1,\dots, b_m$ of
$\frak g_0$ we have $r = \sum r_{ij}b_i \wedge b_j$ then
$F(b_i,b_j) = (r^{-1})_{ij}$ is a non-degenerate
2-cocycle. Conversely, if $\frak f$ is quasi-Frobenius with non-degenerate 2-cocycle $F$ then the inverse of
$F$ is a solution $r$ to the CYBE whose carrier is $\frak f.$
Finally, $r$ and $r'$ with corresponding $(\frak f, F)$ and
$(\frak f',F')$
are equivalent if there is an inner automorphism of
$\frak g$ carrying $\frak f$ to $\frak f'$ and carrying $F$ to a
2-cocycle of $\frak f'$ cohomologous to $F'.$  If
$\frak f$ equals $\frak f'$ and is Frobenius then $r$ and $r'$
are equivalent, so solutions to the CYBE with Frobenius carriers
can be classified by their carriers, but the problem remains
difficult. The classification of quasi-Frobenius subalgebras is
harder yet since, for example, all even-dimensional abelian subalgebras are  
quasi-Frobenius.
 Stolin has carefully used the preceeding analysis
to list (up to equivalence) all of the quasi-Frobenius
subalgebras of $\frak {sl}(2)$ and $\frak {sl}(3)$, see [S2].
He gives the unique solution, $(e_{11}-e_{22})\wedge e_{12}$, of the CYBE
which lies in $\frak {sl}(2)\bigwedge \frak {sl}(2)$, but he
does not explicitly exhibit any which lie in
$\frak {sl}(3)\bigwedge \frak {sl}(3).$

In view of this we consider a restricted problem. If $r$
is a solution to the MCYBE and we identify it with $\lambda \cdot r$
for any scalar $\lambda$, then we can view this equivalence class
as a point inside
the projective space
$\Bbb P(\frak g \bigwedge \frak g)$ of lines in $\frak g \bigwedge \frak
g$. We will show
in Theorem 4.1 that the set of solutions, $\Cal M$, to the MCYBE
is a quasi-projective variety, in other words $\Cal M$ is an
open subset of its closure. The boundary points of $\Cal M$
must be solutions to the CYBE. We call these {\it
boundary solutions,} and ask if it is possible to classify
constructively at least these. We show in Theorem 4.3
that the boundary solutions form a proper subset of the solutions  
to the CYBE
by explicitly constructing non-boundary solutions for
for the classical simple Lie
algebras.

In section 5, we present some families of boundary solutions
to the CYBE lying in
$\frak {sl}(n)\bigwedge \frak{sl}(n).$ The most intriguing of these
are associated to maximal parabolic subalgebras of $\frak {sl}(n)$.
To every $i$ between $1$ and $n-1$ we
can associate a maximal parabolic subalgebra $\frak p_i$ of $\frak  
{sl}(n)$,
namely, that generated by the Cartan subalgebra, all
positive root vectors, and all negative root vectors except
$e_{i+1,i}$. It is known that $\frak p_i$ is Frobenius if and only
if $i$ and $n$ are relatively prime, see [E1]. Also $H^2(\frak p_i)=0$,
see [F], and so a maximal parabolic subalgebra is
Frobenius if and only if it is quasi-Frobenius. Thus $\frak p_i$ is the
carrier of a solution to the CYBE if and only if $(i,n)=1$ in which
case there is a unique solution (up to equivalence) in $\frak p_i\bigwedge \frak p_i$.
We conjecture that in this case the solution is in fact a boundary
solution and lies in the closure of a suitable orbit of the
associated generalized Cremmer-Gervais solution to the MCYBE;
we show that this is indeed true for the ``end''
(maximal) parabolic subalgebras $\frak p_{1}$ and
$\frak p_{n-1}$, see Theorem 5.9 where we give the formula for
this boundary solution. When $n=3$, our formula gives
$$\left(\tsize {\frac{2}{3}}e_{11}-\tsize {\frac{1}{3}}e_{22}-
\tsize {\frac{1}{3}}e_{33}
\right)\wedge e_{12}+ \left( \tsize {\frac{1}{3}}e_{11}+\tsize  
{\frac{1}{3}}e_{22}-\tsize {\frac{2}{3}}e_{33}\right)
\wedge e_{23}+e_{13}\wedge e_{32}\tag 0.2$$
which clearly has carrier $\frak p_1$, the subalgebra
of $\frak {sl}(3)$ consisting of traceless matrices of the form
$$\pmatrix
*&*&*\\
0&*&*\\
0&*&*\endpmatrix.$$
It seems remarkable that this solution to the CYBE can be constructed
directly from (0.1), the associated Cremmer-Gervais solution to the
MCYBE. We also show that the Cremmer-Gervais solution to the MCYBE and the
boundary solution of the CYBE with carrier $\frak p_1$ are part of a
natural simple three-dimensional submodule of
$\frak {sl}(n)\bigwedge \frak {sl}(n)$ under the adjoint action
of the principal three-dimensional subalgebra of $\frak {sl}(n)$.

We can prove our conjecture also for the case
$n=5, i=2$ and show how these results might be extended to establish
the general case, see Conjecture 6.1. 
Finally, we close the paper with a brief discussion about
constructing solutions to the {\it {quantum}} Yang-Baxter equation from
solutions to the CYBE.

Our main interest in solutions to the  MCYBE and CYBE is their connection
to Poisson geometry and quantizations.
If
$r\in \frak g_0 \bigwedge \frak g_0$ is a solution to
the CYBE with carrier $\frak g_0$, then there is
a left invariant symplectic structure on the
corresponding simply connected algebraic group $G_0$. This means that
the coordinate ring $\Cal O(G_0)$ of $G_0$ is equipped with a
left invariant skew bracket,
$\{\hskip 5pt ,\hskip 5pt \}_{r}$ (the ``Poisson bracket''),
which satisfies the Jacobi identity
and has $\{f, gh\}_r=\{f,g\}_rh +\{f,h\}_rg$ for all $f,g,h \in \Cal O(G)$.
If $r=\sum \lambda_{ij}a_i \wedge b_j$ then the associated Poisson  
bracket is
given by 
$\{f,g\}_{r}=\sum \lambda_{ij}(A_i(f)\cdot B_j(g)- B_j(f)\cdot
A_i(g))$ where
$A_i$ and $B_j$ are the left-invariant vector fields on $G$ corresponding
to $a_i$ and $b_j$. 
Similarly, there is a right invariant symplectic structure,
$\{\hskip 5pt ,\hskip 5pt \}'_{r}$
on $\Cal O(G)$ given by
$\{f,g\}_{r}=\sum \lambda_{ij}(A'_i(f)\cdot B'_j(g)-B_j'(f)\cdot
A_i'(g))$ where
$A_i'$ and $B_j'$ are the right-invariant vector fields on $G$  
corresponding
to $a_i$ and $b_j$.

If $G$ is an algebraic group, then it is natural to consider Poisson-Lie
structures on the Hopf algebra
$\Cal O(G)$. A Poisson Lie structure is a Poisson bracket
on $\Cal O(G)$ which is compatible with the comultiplication
$\Delta: \Cal O(G)\to \Cal O(G)\otimes \Cal O(G)$. This
means that $\Delta (\{f,g \})=\{\Delta f, \Delta g\}$ for
all $f,g\in \Cal O(G)$. Drinfel'd proves in [D1] that
all Poisson-Lie structures on $\Cal O(G)$ are of the form
$\{f,g \}= \{f,g\}_r -\{f,g\}'_r$ where $r\in \frak g\bigwedge
\frak g$ satisfies either the MCYBE or CYBE. This bracket,
however, is neither left nor right invariant.

Although we will not discuss deformations in this paper, we would
like to point out that Poisson (and Poisson-Lie) structures correspond
to first order or infinitesimal deformations, which by the above
remarks in turn correspond
to solutions $r$ of the CYBE and the MCYBE. Starting
with such an $r$, one would like, optimally, to have a
way of constructing a deformation with $r$ as its classical limit.
For details of the existence of such deformations, see [BFGP], [D2]
[EK1],and [EK2]. For some recent results on how these deformations
might be constructed see [GGS3] and [GZ].

\subhead 1. The classical Yang-Baxter equations
 \endsubhead
Suppose that $\frak g$ is a finite dimensional simple Lie algebra,
over $\Bbb C$ or split over $\Bbb R$. Throughout this
paper, we will simply use $k$ to denote the ground field. All tensor
products will be taken over $k$, and if $W$ is a $k$-vector
space then we will view $W\bigwedge W$ and $W\bigwedge W\bigwedge W$
as subspaces of $W\otimes W$ and $W\otimes W\otimes W$ in the  
natural way.
In particular, $w_1\wedge w_2$ corresponds to $\frac{1}{2}(w_1\otimes
w_2 - w_2\otimes w_1)$. Let
$U\frak g$ be the
universal enveloping algebra of $\frak g$ and let
$G$ be a simply connected algebraic
group with Lie$(G)=\frak g$.
For $r=\sum a_i\otimes b_i \, \in \frak g \otimes \frak g$
define $\langle r, r \rangle \in U\frak g \otimes U\frak g
\otimes U\frak g$ to be
$$[r_{12},r_{13}]+[r_{12},r_{23}]+[r_{13},r_{23}]$$
where $r_{12}=r\otimes 1$, $r_{23}=1\otimes r$,
$r_{13}=\sum a_i\otimes 1\otimes b_i$.
\definition{Definitions 1.1}
\roster
\item An element $r\in \frak g \otimes \frak g$ is a solution to
the classical Yang-Baxter equation (CYBE) if $\r =0$.
\item  An element $r\in \frak g \otimes \frak g$ is a solution to
the modified classical Yang-Baxter equation (MCYBE) if $\r$ is non-zero
and $\frak g$-invariant.
\item A solution to either the CYBE or MCYBE is unitary if
$r$ is skew-symmetric and non-unitary otherwise.
\endroster
\enddefinition
Actually, what we have just defined are constant solutions to the
Yang-Baxter equations. There are related non-constant
solutions to these equations, (those which depend on a spectral  
parameter),
but our only concern here is with the constant solutions. We will also
only be considering unitary solutions to the CYBE and MCYBE.
As discussed in the introduction, it is these solutions which have  
important
homological and geometrical meanings.
A natural question therefore is to determine
and describe, if possible, the set of all $r\in  \frak g \bigwedge  
\frak g$
which satisfy either the CYBE or the MCYBE. In our analysis of this  
question,
it will be useful to identify $r \in \frak g \bigwedge \frak g$ with
$\lambda \cdot r$ for any scalar $\lambda \in k ^\times$ and so if
$d={\text {dim}}_{k}\,(\frak g)$, we can view $r$ as an element of the
${{d\choose 2}-1}$ dimensional projective space
$\Bbb P(\frak g \bigwedge \frak g)$. Similarly we can view $\r$, which
necessarily lies in $\frak g \bigwedge \frak g \bigwedge \frak g$, as
an element of the ${{d\choose 3}-1}$ dimensional projective
$\Bbb P(\frak g \bigwedge \frak g \bigwedge \frak g)$.
\definition{Notation} Let $\Cal C$ and $\Cal M $ denote, respectively,
the subsets of $\Bbb P(\frak g \bigwedge \frak g)$
consisting of solutions to the CYBE and MCYBE.
\enddefinition
There
is an important notion of equivalence of solutions to these
equations which comes from the
diagonal action of the group $G$ on $\frak g \bigwedge \frak g$.
Specifically, if $g\in G$ and $r=\sum a_i\wedge b_i$ then we say
that $r$ and $g\cdot r =\sum ga_ig^{-1}\wedge gb_ig^{-1}$ are equivalent.
This notion of equivalence is well-defined
since both $\Cal M$ and $\Cal C$ are invariant under
under the action of $G$.
\subhead{2. Solutions of the modified classical Yang-Baxter equation}
\endsubhead
We first discuss $\Cal M$,
the space of solutions to the MCYBE. Remarkably, Belavin and Drinfel'd  
have given a constructive description of $\Cal M$ in [BD].
They show, in particular, that $\Cal M$ is a finite disjoint union of
components each of which is determined by certain data associated with
the root system of $\frak g$. To describe their work, we first need to
recall some important facts and notation pertaining to finite
dimensional simple Lie algebras.
Let $(\hskip 4pt ,\hskip 4pt)$ be an invariant non-degenerate symmetric
bilinear form on $\frak g$, let $\frak h \subset \frak g$
be a fixed a Cartan subalgebra, and
denote the root system by $\Phi$.
Write $\Phi =\Phi_+ \cup \Phi_-$ where
$\Phi_+$ and $\Phi_-$ are the positive and negative roots,
respectively, and let $\Pi \subset \Phi_+$
be a basis of positive simple
roots.
There is then a basis
$\{ x_{\tau}\, | \,\tau \in \Phi _+\} \cup \{x_{-\tau}\,|\, \tau  
\in \Phi _+\}
\cup \{h_{\tau} \,| \, \tau \in \Pi \}$ of $\frak g$ for which
$(x_{\tau}, x_{-\tau})=1$ for all $\tau \in \Phi_+$.
This gives a triangular decomposition
$\frak g =\frak n_- \oplus \frak h \oplus \frak n_+$ where $\frak  
n_-$ and $\frak n_+$ are the nilpotent
subalgebras spanned by the negative and positive root vectors
relative to the chosen Cartan subalgebra.
The data which determine the components of $\Cal M$ are called,
using the terminology of [BD],
admissible triples.
\definition{Definition 2.1 [BD]} An {\it {admissible triple}} is a
triple $\Cal T=(\Pi_1,\Pi_2,T)$
where $\Pi_1$
and $\Pi_2$ are subsets of $\Pi$, with
$\# (\Pi_1)=\# (\Pi_2)$
and $T:\Pi_1\rightarrow \Pi_2$ is a bijection with the properties that
\roster
 \item $T$ preserves the Killing form, {\it i.e.,}
 $(T(x_\pi),T(x_\rho)) = (x_\pi,x_\rho)$ for all $\pi,\rho \in  
\Pi_1$, and
 \item for every $\pi \in \Pi_1$ there is a positive integer $m$
 such that $T^m\pi \notin \Pi_1.$
 \endroster \enddefinition
The following theorem proved in [BD] (which we state in a manner
convenient for our purposes), shows that
the set of admissible triples gives a decomposition of $\Cal M$
into mutually disjoint
components.
\proclaim{Theorem 2.2 ([BD])} Let
$\Cal M\in \Bbb P(\frak g\bigwedge \frak g)$ be the set of  
solutions to the
MCYBE.
 For every admissible triple $\Cal T$, there is a
non-empty subset $\Cal M_{\Cal T}$ of $\Cal M$ such that the following hold:
\roster
\item $\Cal M = \bigcup \Cal M_{\Cal T}$ where $\Cal T$
runs through the set of all
admissible triples.
\item If $\Cal T\neq \Cal T'$, then
$\Cal M_{\Cal T}\bigcap \Cal M_{\Cal T'}=\emptyset$.
\endroster
\endproclaim
To each triple $\Cal T$ there is a certain
affine
subvariety $\Cal B_{\Cal T}\subset \frak h \bigwedge \frak h$ which
determines the dimension of $\Cal M_{\Cal T}$. Its description,
given in the next theorem, uses for every 
$\pi \in \Pi_1$ the map $1\otimes \pi :\frak h \otimes
\frak h \rightarrow \frak h$ defined by
$(1\otimes \pi)(h_\tau \otimes h_\rho)=\pi (h_\rho)\, h_\tau.$
In addition to $\Cal B_\Cal T$, there is a uniquely defined
element $\alpha \in \frak n_+ \bigwedge \frak n_-$ 
determined by $\Cal T$ which is used to describe $\Cal M_\Cal T$. Its
definition uses
the linear extension of $T$ to the set
$\widehat{\Pi}_1 \subset \Phi_+$ of positive roots which are sums  
of simple
roots in $\Pi_1$. If $\pi \in \widehat\Pi_1$ and $\varrho \in \Phi_+$, then
define $\pi \prec \varrho$ if there is an $m > 0$
such that $T^m\pi = \varrho.$
\proclaim{Theorem 2.3 ([BD])} Let $\Cal T=(\Pi_1, \Pi_2,T)$ be an
admissible triple and let $\Cal B_\Cal T$ be the set of all
$\beta \in \frak h \bigwedge \frak h$ for which
$$(1\otimes (T(\pi )-\pi))\, \beta = {\tsize {\frac{1}{2}}}(h_{T(\pi)}+
h_\pi) \qquad {\text {for all}} \qquad \pi \in \Pi_1.$$
\roster
\item  $\Cal B_\Cal T $ is an affine variety of dimension
$\binom {d}{2}$ where $d=\#(\Pi_1 -\Pi).$ 
\item If
$$\gamma =\sum _{\pi \in \Phi _+}x_{\pi}\wedge x_{-\pi} \quad
{\text {and}} \quad
\alpha = 2\sum _{\Sb \pi \prec \varrho\\
\pi,\varrho \in \Phi_+\endSb}x_{\pi} \wedge x_{-\varrho},$$
then $\gamma +\beta +\alpha$ is a solution to the MCYBE
for any $\beta \in \Cal B_\Cal T$.
\item Every solution of the MCYBE is equivalent to a unique
solution of the form $\gamma +\beta +\alpha$ associated to some triple
$\Cal T$.
\endroster
\endproclaim
The previous two theorems thus give a complete
constructive description of $\Cal M$ up to equivalence; the triple
$\Cal T$ is the ``discrete'' parameter and the variety $\Cal B_\Cal T$
is the ``continuous'' parameter of of the component $\Cal M_\Cal T$.
These theorems do not, however, give a direct way of finding
specific $\beta \in \Cal B_\Cal T$. For some triples associated
to ${\frak {sl}}(n)$ we are able to completely describe $\Cal B_\Cal T$,
see Theorem 2.9.
Details of our remarks about the
Belavin-Drinfel'd classification can be found in either [BD] or
[CP, \S 3.1.A and \S 3.1.B].

Before going on to discuss $\Cal C$, we will examine
some important triples and their associated solutions to the MCYBE.
First, let us recall the natural
partial ordering of triples discussed in [GGS3, \S 6].
If $\Cal T=(\Pi_1,\Pi_2,T)$ and
$\Cal T'=(\Pi_1',\Pi_2',T')$ are triples then set
$\Cal T < \Cal T'$ if $T$ is the restriction of $T'$ to
some subset $\Pi_1$ of $\Pi_1'$.
The trivial triple, $\Cal T_{\emptyset}$, is the smallest element
in this ordering and has $\Pi_1 =\Pi_2 =\emptyset$.
The $\alpha$-term of the corresponding solution
is zero while the
$\beta$-term can be arbitrarily chosen from $\frak h \bigwedge  
\frak h$ and
so if $l$ is the
rank of $\frak g$, this determines, up to equivalence,
 an ${l \choose 2}+1$ dimensional family of solutions to the MCYBE.
The dimension
of this family of solutions is greater than the dimension of any
family of solutions associated to a non-trivial triple.

The simplest example of a non-trivial triple occurs when
$\frak g= \frak {sl}(3)$.
The two simple roots of this Lie algebra may be identified
 with the set $\{1,2\}$ and the associated triple has
$\Pi_1=\{1\}$, $\Pi_2=\{2\}$, and $T(1)=2$.
  Here, $\beta $ is unique since $\Pi _1$
omitted all but one positive simple root. The corresponding
solution $\gamma +\beta +\alpha $ is the Cremmer-Gervais
solution to the MCYBE, see (0.1).
This triple has an analog in $\frak{sl}(n)$ for
 all $n>3$. If we make the natural
identification of $\Pi$ with $\{1,2,\ldots n-1\}$ and
set $\Pi_1=\{1,2,\ldots n-2\}$ and $\Pi_2=\{2,3,\ldots n-1\}$, then
the triple
$\Cal T_{CG}$ (notation explained below)
determined by the rule $T(i)=i+1$ is admissible. This triple
is maximal in the given partial ordering because $\Pi_1$ omitted
only one root. Morever,
the $\beta$-term for the associated solution to the MCYBE is
uniquely determined and hence there is a unique solution
$r_{CG}$ (given in (5.8)) to the MCYBE associated to this
triple. We call this the {\it {Cremmer-Gervais solution}} because
the associated non-unitary
solution to the CYBE is the classical limit of the
quantum
Yang-Baxter matrix exhibited in [CG].
In contrast with the standard
quantum groups and their multiparameter versions, little is known
about the Cremmer-Gervais quantum groups. See [H] for some
recent progress in understanding
their structure.

For $\frak {sl}(n)$, there are other
triples for which, like the Cremmer-Gervais
triple, $\Pi_1$ consists of all but one positive simple root of  
$\frak {sl}(n)$.
For each of these, the $\beta$-term is uniquely
determined and so, up to equivalence, they each correspond to a unique
solution to the MCYBE. We call these the {\it {generalized  
Cremmer-Gervais}}
triples. These are clearly maximal elements in the partial ordering
of the triples mentioned earlier.
If we again make the identification of $\Pi$ with the set
$\{1,2,\ldots ,n-1\}$ then the conditions for building
an admissible generalized Cremmer-Gervais triple
translate into finding a bijection $T:S_1\rightarrow S_2$ between
subsets $S_1$ and $S_2$ of $\{1,2,\ldots ,n-1\}$ such that:
\roster
\item "(2.4)"For every $i \in S_1$
there is an $r$ such that $T^ri \notin S_1$, and
\item "(2.5)"$T$
respects adjacency, i.e., if $i,j \in S_1$ then $|i-j| = 1$
implies $|Ti - Tj| = 1,$ and conversely.
\endroster
As the next theorem shows, the number of 
generalized Cremmer-Gervais triples depends
on the factorization of $n$.

\proclaim{Theorem 2.6} Suppose that $S_1$ and $S_2$ are
subsets of $\{1, \dots, n-
1\}$ with $\#S_1 = n-2$ and that $T:S_1 \to S_2$ determines a
generalized Cremmer-Gervais triple.
Let the omitted element of $S_1$ be $n-i.$ Then $i$ and $n$ are
relatively prime,
the omitted
element of $S_2$ is $i$, and $T$ sends every $j \in S_1$ to $j+i
\mod n.$
\endproclaim

\remark{Remark 2.7} A consequence of the Theorem is that, modulo $n$,
$S_1=\{i,\, 2i,\, \ldots ,\,(n-2)i\}$, $S_2=\{2i,\, 3i\,
\ldots ,\, (n-1)i\}$ and these sets have a natural
order (this order, however, is
not compatible with the order in $\Bbb N$).  
\endremark
\demo{Proof of 2.6} Consider first the case where $i=1,$ so $S_1 =
\{1,\dots,n-2\}.$ The matter is trivial for $n=3$ so we may
suppose $n > 3.$ Since $T$ preserves adjacency, we must have $S_2
= \{2,\dots,n-1\},$ else $T$ would be a bijection of $S_1$ with
itself and condition 2.4 would be violated. The
only question then is whether $T$ preserves order, i.e., sends every
$j$ to $j+1$ (as we claim), or reverses order, sending $j$ to $n-
j.$  The latter is impossible, for $\{2,\dots,n-2\}$ would then
be a non-empty subsegment of $S_1$ carried onto itself by $T,$
and this would violate condition 2.4. Exactly the
same argument applies when $i = n-1,$ so we may now suppose that
$1 < i < n-1.$  Removing $n-i$ then breaks the string of integers
$1,\dots,n-1$ into two non-empty segments, $S_1' = \{ 1,\ldots ,  
n-i-1 \} $
and $S_1'' = \{n-i+1,\ldots n-i\}.$  Their images must again be separated
segments of $\{1,\dots,n-1\},$ for otherwise some pair of
non-adjacent elements of $S_1$ would be carried by $T$ to adjacent
elements and condition 2.5 would be violated.
Now the only way that $S_1'$ and $S_1''$ can be sent to
non-adjacent and non-overlapping segments of $\{1,\dots, n-1\}$
without either being sent into itself is if $S_1'$ is sent to a
terminal segment and $S_1''$ to an initial segment of
$\{ 1,\dots,n-1 \}$ i.e., that the image of $S_1'$ is
$\{i+1,\dots,n-1\}$ and that of $S_1''$ is $\{1,\dots,i-1\}.$  In
particular, the element omitted from $S_2$ must be $i.$
Moreover, we can not have $i=n-i$ or $T$ would carry $S_1$ onto
itself.

Now the possibilities for $T$ are that it either preserves or
reverses the order in $S_1',$ i.e., that if $j \in S_1'$ then
either $T(j) = j+i$ or $T(j) = n-j,$ and similarly for $S_1''.$  One
sees now that $i$ and $n$ must be relatively prime.
For suppose that $(n,i)= m > 1.$  Then in each of the four
possible cases, $T$ carries a
multiple of $m$ lying in $S_1$ to another multiple of $m,$ and
the complement of the set of multiples of $m$ is not empty. This
complement is contained in $S_1,$ so we would have a subset of
$S_1$ carried bijectively onto itself by $T,$ which is
impossible.  All that remains, therefore, is to show that $T$
preserves the order separately in $S_1'$ and in $S_1''.$  Note
that these segments are non-empty and of different lengths.  We
must show that it is impossible for $T$ to reverse the order
either in $S_1'$ or in $S_1''$ or in both.  The last case is
trivial, for then $T1 = n-1$ and $T(n-1) = 1,$ so condition 2.4
would be violated.  The arguments for the other cases are similar,
so we shall do only one explicitly. Suppose that $S_1' =
\{1,\dots,n-i-1\}$ is the longer segment, and that $T$ reverses
its order, sending any $j$ with $1 \leq j \leq n-i-1$ to $n-j.$
It follows that $i+1 \leq n-i-1$ so $\{i+1,\dots,n-i-1\}$ is a
non-empty subsegment of $S_1'$ which is carried by $T$ onto
itself, which is impossible. \qed
\enddemo
Let $\Cal T_i$
denote the generalized Cremmer-Gervais
triple described by Theorem 2.6. As stated earlier,
$\Cal M_{\Cal T_{i}}$ contains, up to equivalence, a unique solution to
the MCYBE because there is only a single $\beta \in \frak h \bigwedge \frak h$ 
which satisfies 
$$(1\otimes (T(\pi )-\pi))\, \beta = {\tsize {\frac{1}{2}}}(h_{T(\pi)}+
h_\pi) \qquad {\text {for all}} \qquad \pi \in \Pi_1.\tag 2.8$$
The following explicitly determines this element. 
Since $\frak g ={\frak {sl}}(n)$ we
make the identifications of $\Pi$ with $\{1,2,\ldots n-1\}$,
of $x_j$ with $e_{j,j+1}$, and of $h_j$ with $e_{jj}-e_{j+1,j+1}$.
\proclaim{Theorem 2.9} Suppose that $i$ and $n$ are
relatively prime and that $\Cal T_i$ is the generalized Cremmer-Gervais
triple.
Let $\beta = \sum _{l<j}b_{p,q}e_{pp}\wedge e_{qq}=(1/2)\sum _{p<q}
b_{p,q}(e_{pp}\otimes e_{qq}-e_{qq}\otimes e_{pp})$ be the unique solution
to (2.8). If $p-q=si \mod n$ then $b_{p,q}=(1/n)(n-2s)$.
\endproclaim
\demo{Proof}
Recall that for $\Cal T_i$ we have
$\Pi_1=\{i,\, 2i,\, \ldots ,\,(n-2)i\},\quad \Pi_2=\{2i,\, 3i\,
\ldots ,\, (n-1)i\},$ and $T:\Pi_1\rightarrow \Pi_2$ is given by
$ T(ri)=(r+1)i \mod n$ (see Theorem 2.6).
With the indicated identifications above, finding a $\beta \in \frak h
\bigwedge \frak h$ which satisfies
(2.8) means we must have
$$(1\otimes (e_{(r+1)i,(r+1)i+1}-e_{ri,ri+1}))\beta = 
{\tsize {\frac {1}{2}}}(e_{(r+1)i,(r+1)i}-e_{(r+1)i+1,(r+1)i+1}+
e_{ri,ri}-e_{ri+1,ri+1})\tag 2.10$$
for all $ ri\in \Pi_1.$ For $p<q$, set $b_{q,p}=-b_{p,q}$. Then
an elementary computation shows that
$$(1\otimes e_{si,si+1})\, \beta =
\left( \sum _{j\neq si,si+1}(b_{j,si}-b_{j,si+1})e_{jj}\right)
-b_{si,si+1}(e_{si,si}+e_{si+1,si+1}).\tag 2.11 $$
The result now follows by combining (2.10) and (2.11) 
with the given formula for the coefficients
$b_{p,q}$.  \qed \enddemo
Theorems 2.3, 2.6, and 2.9 thus provide
an explicit
description of all generalized Cremmer-Gervais triples
and their corresponding solutions to the MCYBE. 

\subhead{3. Solutions of the classical Yang-Baxter equation and
quasi-Frobenius Lie algebras}
\endsubhead
The space $\Cal C$ of skew solutions to the CYBE
is, in a sense, less
understood than $\Cal M$. In fact, as Belavin and Drinfel'd remark,
a constructive classification of the solutions to the CYBE is essentially
intractable as it would, in particular, require knowledge
of all abelian subalgebras $\frak a \subset \frak g$ since any
$r\in \frak a  \bigwedge \frak a$ has $\r =0$. On the positive side,
there is a
homological interpretation of
the unitary solutions to the CYBE, see [BD, Proposition 2.4].
In [S1] and [S2], Stolin
used this connection
to provide a non-constructive description of $\Cal C$ in terms
of {\it {quasi-Frobenius}} Lie algebras.
\definition{Definition 3.1}\roster
\item A Lie algebra $\frak f$ is
{\it Frobenius} if there exists a linear functional $f\in \frak f^*$
such that the skew bilinear form
$[ \hskip 4pt, \hskip 4pt]_f:  \frak f\bigwedge \frak f \rightarrow k$
which sends $x\wedge y$ to $f([x,y])$
is non-degenerate (here $k$ is $\Bbb C$ or $\Bbb R$). Equivalently, 
$[ \hskip 4pt, \hskip 4pt]_f$ is a non-degererate two-coboundary
 for the Lie algebra cohomology of
$\frak f$ with coefficients in $k$.
\item A Lie algebra $\frak f$ is {\it{quasi-Frobenius}} if there exists
a non-degenerate two-cocycle $F:\frak f\bigwedge \frak f \rightarrow k$.
(Note that a Frobenius Lie algebra is also quasi-Frobenius.)
\endroster \enddefinition
\remark{Remark 3.2} A result of A. I. Ooms states that any finite
dimensional Frobenius Lie algebra
has a primitive enveloping algebra. The converse holds in case the Lie
algebra is algebraic, see [O].
\endremark

The connection between quasi-Frobenius Lie algebras and the CYBE is  
that if
$F_{ij}$ is the matrix of the form $F$ relative to some basis
$x_1,\ldots x_d$ of $\frak f$ then
$r=\sum (F^{-1})_{ij}x_i\wedge x_j\in \frak f\bigwedge \frak f$
is a skew solution to the CYBE.
Conversely, if we start with a non-degenerate
skew solution $\sum r_{ij}x_i\wedge x_j\in \frak f\bigwedge \frak f$ 
to the CYBE, then the map $F: \frak f \bigwedge \frak f\rightarrow k$
where
$F(x_i,x_j)=(r^{-1})_{ij}$ is a non-degenerate two-cocycle and
so $\frak f$ is quasi-Frobenius.
A simple Lie algebra $\frak g$
can not itself be quasi-Frobenius and so any
solution $r\in \frak g \bigwedge \frak g$ is necessarily degenerate.
For any such $r$, however there always exists a largest subalgebra
$\frak f\subset \frak g$ on which $r$ is non-degenerate,
see [CP, Propositions 2.2.5 and 2.2.6]. That
subalgebra is necessarily quasi-Frobenius. Thus, finding solutions
to the CYBE in $\frak g\bigwedge \frak g$
is equivalent to finding pairs $(\frak f, F)$ where $\frak f$ is a  
subalgebra
of $\frak g$ and $F$ is a non-degenerate two-cocycle on $\frak f$.
We call $\frak f$ the
{{\it carrier}} of $r$.
If $(\frak f, F)$ and $(\frak f', F')$
are quasi-Frobenius subalgebras and
$\phi :\frak g\rightarrow \frak g$ is an inner
automorphism which takes $\frak f$ to $\frak f'$,
then the solutions to the CYBE associated to $(\frak f,F)$
and $(\frak f',F')$ are equivalent if and only if
$\phi^{*}F'-F$ is cohomologous to zero. In particular, if
$\frak f=\frak f'$ and is Frobenius, then up to equivalence
it contains a unique solution to the CYBE. It seems difficult
though to classify all Frobenius subalgebras of
simple Lie algebra $\frak g$. The same problem for quasi-Frobenius
subalgebras is even harder, for example,
any even-dimensional abelian subalgebra is
quasi-Frobenius. Some examples of Frobenius Lie algebras and a discussion
of quasi-Frobenius Lie algebras can be found in [E2].
\subhead{4. Boundary solutions to the classical Yang-Baxter equation}
\endsubhead
Even though the Belavin-Drinfel'd classification of solutions
to the MCYBE and the homological
interpretation of solutions of the CYBE have been known for some  
time now,
there has been little attention paid to the interface between these
equations. As we shall see next, certain solutions to the CYBE
can be viewed as limiting cases of solutions to the MCYBE.
 \proclaim{Theorem 4.1} Let $\Cal C$ and $\Cal M$ be the subsets
of $\Bbb P(\frak g\bigwedge \frak g)$ consisting of solutions to the
CYBE and MCYBE, respectively, and denote their Zariski closures
by $\overline {\Cal C}$ and $\overline {\Cal M}.$
\roster
\item $\Cal C =\overline {\Cal C}$, i.e. $\Cal C$ is a (closed) variety.
\item  $\Cal M$ is a quasi-projective
variety, i.e. it is an open subset of $\overline {\Cal M}$.
\item  Any point on the boundary of $\Cal M$ lies in $\Cal C$
and hence is a solution to the CYBE.
\endroster \endproclaim
\demo{Proof}
(1) Suppose $r=\sum _{p,q}r_{pq}\, x_p\wedge x_q$ and $\r=
\sum_{i,j,k}c_{ijk}\, x_i\wedge x_j\wedge x_k$. It is easy to see
that each $c_{ijk}=f_{ijk}(r)$ for some homogeneous quadratic
polynomial function $f_{ijk}:\Bbb P(\frak g\bigwedge \frak  
g)\rightarrow k$
and so $\Cal C$ is a (closed) variety.
\newline (2) and (3) Consider the map which sends $r$ to
$\r$; it is a well-defined map $\phi : (\Bbb P(\frak g \bigwedge \frak g)
-\Cal C)
\rightarrow \Bbb P(\frak g \bigwedge \frak g \bigwedge \frak g)$.
Since $\frak g$ is simple, the space of invariants in
$\frak g \bigwedge \frak g \bigwedge \frak g$ consists
of all non-zero multiples of some fixed non-zero invariant in
$\frak g \bigwedge \frak g \bigwedge \frak g$.
As a subspace of
$\Bbb P(\frak g \bigwedge \frak g \bigwedge \frak g)$, this is just
a single point which we denote $\omega$.
Therefore
$\Cal M =\phi ^{-1}(\omega)$ and so
$\Cal M$ is a closed subset of $(\Bbb P(\frak g \bigwedge \frak g)
-\Cal C)$. Hence every element
of $\overline {\Cal M} - \Cal M$ must lie in $\Cal C$. Finally, to
prove that $\Cal M$ is open in $\overline {\Cal M}$, it suffices to 
show that $\overline {\Cal M} - \Cal M$ is a closed
subset of $\overline {\Cal M}$. But this is trivial since
$(\overline {\Cal M} - \Cal M)=\Cal C \cap \overline {\Cal M}$ and
both $\Cal C$ and $\overline {\Cal M}$ are closed. \qed
\enddemo
In light of the preceeding theorem, we can focus on a more restricted
class of solutions to the CYBE.
\definition{Definition 4.2} Let $\Cal C$ and $\overline {\Cal M}$
be defined as in Theorem 4.1. An element $r\in \frak g\bigwedge \frak g$
is a {\it {boundary solution}} to the CYBE if
$r\in\Cal C\cap  \overline {\Cal M}.$, i.e. if $r$ is in the closure
of, but not contained in the space of solutions of the MCYBE.
\enddefinition
It is not a priori clear that there exist non-boundary solutions, but
using the Belavin-Drinfel'd classification of solutions to the MCYBE, we
can at least show that
for the classical simple Lie algebras,
there do indeed exist non-boundary solutions
of the CYBE.
\proclaim{Proposition 4.3} Let
$\frak g$ be a classical simple Lie algebra and
set $W=\frak g \bigcap X$ where $X$ is the space
spanned by all matrix units $e_{ij}$ with
$ i\leq [{n\over 2}]$ and
$j\geq [{n\over 2}]$. Then, if the rank of $\frak g$ is large,
the generic
element of $W\bigwedge W$ is a non-boundary solution of the CYBE.
\endproclaim
\demo{Proof} We only give the proof for $\frak g=\frak{sl}(n)$ since it
easily modifies to the other cases. Since $W=X$ for $\frak  
g=\frak{sl}(n)$,
the dimension of $W\bigwedge W$ is on
the order of $n^4/32$ for large $n$. Moreover,
since $W$ is abelian, every element of $W\bigwedge W$
satisfies the CYBE. Now the component of the solutions to the MCYBE
of maximal dimension corresponds to the triple with
$\Pi_1 =\Pi_2 =\emptyset$. The associated solution to the MCYBE
has $\alpha =0$ and
$\beta \in \frak h \bigwedge \frak h$ can be arbitrary. Now since
${\text {dim}}({\frak {sl}}(n))=n^2-1$ and ${\text {dim}}(\frak h  
\bigwedge
\frak h)={{n-1}\choose 2}$, this component has dimension at most
$(n^2-1)+{1\over 2}(n-1)(n-2)$ which is on the order of ${3\over 2} n^2$
for large $n$. Hence the generic element of $W\bigwedge W$ must be
a non-boundary solution of the CYBE
since its dimension exceeds ${3\over 2} n^2$.
\qed \enddemo

\subhead{5. Construction of boundary solutions}
 \endsubhead
Although we presently can not classify, or even find all boundary  
solutions,
we are able to exhibit
several (previously unknown) families of solutions
to the CYBE.
We will do this
by analyzing the action of the group $G$ on certain solutions to the 
MCYBE obtained from Theorem 2.3.
The most
intriguing of these is found using the Cremmer-Gervais solution to
the MCYBE and is
associated to the maximal parabolic subalgebra
of $\frak {sl}(n)$ obtained by deleting the first negative root.
This family, and all others we have, can be found using the following
elementary (but useful) result. In the following, we treat $t$ as a
formal variable and enlarge the coefficient ring to $k[t]$.
\proclaim{Propostion 5.1} Suppose that $r\in \frak g \bigwedge  
\frak g$ and
$r_t=(r+tr_1+\cdots +t^mr_m)\in
\frak g \bigwedge \frak g $ are solutions of the MCYBE
with $\r=\langle r_t,r_t\rangle$.
Then $r_m$ is a boundary solution
to the CYBE.\endproclaim
\demo{Proof}
Since
$\langle r_t,r_t\rangle =\r+t(\langle r_1,r\rangle +\langle  
r,r_1\rangle)+\cdots
+t^{2m}\langle r_m,r_m \rangle$ and $\r=\langle r_t,r_t\rangle$,
the coefficients of each of the terms with positive degree in $t$ must 
vanish identically. In particular, we must have $\langle r_m,r_m  
\rangle=0$
and so $r_m$ satisfies the CYBE.
Now dividing $r_t$ by $t^m$ we obtain
$(r/t^m)+(r_1/t^{m-1})+\cdots +r_m$ which satisfies the MCYBE and,  
any polynomial function which vanishes on $r_t/t^m$ must also vanish  
on
$r_m$ and so $r_m$ is a boundary solution.
\qed \enddemo
An efficient way to construct boundary solutions of this sort
is to consider the orbit of a fixed solution, $r$, of the MCYBE
under the action of $\exp (ta)\in G$ with $a\in \frak g$ nilpotent.
Since $\r$ is invariant it follows that $r$ and $r_t=\exp (t\,a)\cdot r$
satisfy the hypotheses of the proposition.
This idea was first considered in [GGS1, \S 15] and was used to
produce the first family of boundary solutions to the CYBE.
\definition{Example 5.2 (see [GGS1] \S 15)}
Suppose that $\frak g=\frak {sl}(n)$ and
$r=\gamma =\sum _{i<j}e_{ij}\wedge e_{ji}$. Note that $\gamma$ lies  
in the component
$\Cal M_{\Cal T_{\emptyset}}$ of $\Cal M$ corresponding to the
trivial triple. An elementary computation shows that
$$\exp (-t\, e_{1n})\cdot \gamma = \gamma + t\left( (e_{11} - e_{nn})\wedge
 e_{1n} + 2\sum_{i=2}^{n-1}e_{1i}\wedge e_{in}\right)\tag 2.2.1$$
and so by
Proposition 5.1, the coefficient of the linear term,
$(e_{11} - e_{nn})\wedge
 e_{1n} + 2\sum_{i=2}^{n-1}e_{1i}\wedge e_{in}$ is a boundary  
solution to the
CYBE, (in [GGS1] this solution was denoted
$\gamma_{\infty}$, but this notation no
longer seems suitable since, as we shall see, there are many other  
boundary
solutions.) The carrier of this solution
is the Lie algebra, $\Cal H$, consisting of matrices of the form
$$\pmatrix 1&*&\hdots &*&*\\
0&0&\hdots &0&*\\
\vdots&\vdots &\ddots &\vdots &\vdots\\
0&0&\hdots &0&*\\
0&0&\hdots &0&-1\endpmatrix $$
and so is a one-dimensional extension of a Heisenberg algebra.
(Arbitrary scalars can be substituted into those entries marked
with $*$.) Since $\Cal H$ is the carrier Lie algebra of a solution
to the CYBE, it follows that $\Cal H$ is quasi-Frobenius.
In this case though, it is easy to see that $\Cal H$ is actually
Frobenius. The skew bilinear form
$[\hskip 4pt ,\hskip 4pt ]_{e_{1n}^*}$ associated to the
linear functional $e_{1n}^*:\Cal H\rightarrow k$ is non-degenerate.
In fact, the solution to the CYBE obtained from the
inverse of $[ \hskip 4pt ,\hskip 4pt  ]_{e_{1n}^*}$ is the boundary
solution of Example 5.2.
If we take $n=2$ in that example,
we obtain $(e_{11}-e_{22})\wedge e_{12}$,
the unique solution, up to equivalence,
of the CYBE which lies in $\frak {sl}(2)\bigwedge \frak {sl}(2)$.
\enddefinition
We showed in [GGS3, \S 7] that the boundary solution of
Example 5.2 was actually
part of a larger family of solutions to the CYBE. At that
time however, we were unable to determine whether elements of this family
were boundary solutions.
We show next that this
larger family does indeed consist of boundary solutions.
\definition{Example 5.3 (see [GGS3, \S7])}
Let $\frak g =\frak{sl}(n)$ and set $a=\lambda_1 e_{1n}+\lambda_2  
e_{2,n-2}+\cdots +\lambda_d e_{d,n-d+1}$
where $d=[{n\over 2}]$ and $\lambda _i\in k^\times $.
Then
$$\exp (-t\, a)\cdot \gamma =
\gamma +t \left( \sum _{p=1}^{d}\left(\lambda_p^{-1}  
(e_{pp}-e_{n-p+1,n-p+1})
\wedge e_{p,n-p+1}+2\sum_{i=p+1}^{n-p}e_{pi}\wedge e_{i,n-p+1}
\right) \right)$$
and so the coefficient of the linear term is a boundary solution
to the CYBE. The carrier for this solution is the Lie algebra
$\Cal H'$
of matrices of the form
$$\pmatrix d_{1}&*&\hdots &*\\
0&\ddots&\ddots&\vdots \\
\vdots&\ddots &\ddots &*\\
0&\hdots&0 &d_{n}
\endpmatrix $$
where the diagonal entries $d_1,\ldots ,d_n$
are defined as follows:
$d_i=1$ for $i\leq{n\over 2}$, $d_i=-1$ for $i\geq {n\over 2}+1$,
and $d_i=0$ if $n=2i+1$. This Lie algebra is also Frobenius.
The skew form associated to the linear functional
$e_{1n}^*+\cdots e_{d,n-d+1}^*:\Cal H'\rightarrow k$
is non-degenerate and its inverse is the boundary solution
constructed above.
\enddefinition
\remark{Remark 5.4} A slightly more general version
of the preceeding two examples can be obtained by replacing
$\gamma$ with $\gamma +\beta$ with
$\beta \in \frak h \bigwedge \frak h$ arbitrary; recall
that this is the
the generic element of $\Cal M_{\Cal T_{\emptyset}}$.
In a qualitative sense though,
the resulting boundary solutions are essentially identical to those
with $\beta =0$.
The only
difference is that
the carrier Lie algebras
have more general diagonal entries
than just $\pm 1$ or $0$.
\endremark

As stated earlier, the most interesting boundary solutions we have so far
are related to maximal parabolic subalgebras of $\frak {sl}(n)$.
Henceforth we will assume that $\frak g=\frak {sl}(n)$ and
once again we will
make the natural identification of $\{1,\ldots ,n-1\}$
with a fixed basis $\Pi$ of positive simple roots.
To each subset non-empty subset $\Omega \in \Pi$, there
is an associated {\it {parabolic subalgebra}},
$\frak p_{\Omega}$, of $\frak {sl}(n)$
which is generated by the Cartan subalgebra $\frak h$,
all positive simple root
vectors, and
all negative simple root
vectors except those of the form $x_{-\pi}$ with $\pi \in \Omega$.
It is easy
to see that $\frak h \oplus \frak n_{+}\subset \frak p_{\Omega}$
for every $\Omega$. Our primary focus will be on
{\it {maximal}} parabolic subalgebras, those with $\# \Omega =1$.
The maximal parabolic subalgebra corresponding to $\Omega =\{i\}$
will be denoted $\frak p_i$. In this case, the only missing negative
root vector is $e_{i+1,i}$.
The special cases $\frak p_1$ and
$\frak p_{n-1}$ will be called the {\it {end}} (maximal) parabolic
subalgebras, they consist of matrices of the form
$$\frak p_1=\pmatrix *&*&\hdots&*\\
0&*&\hdots &*\\
\vdots&\vdots &\ddots &\vdots\\
0&*&\hdots&*\endpmatrix
\quad {\text{and}}\quad
\frak p_{n-1}=\pmatrix *&*&\hdots &*\\
\vdots &\vdots &\ddots &\vdots\\
*&*&\hdots&*\\
0&\hdots &0&*
\endpmatrix$$
where the only restriction on the entries marked $*$
is that the trace is zero. The end parabolic subalgebras have
have the greatest dimension among all the $\frak p_i$.
Note that $\frak p_i\cong \frak p_{n-i}$
since they map to each other under the automorphism of
$\frak {sl}(n)$ obtained the non-trivial automorphism of its
Dynkin diagram.

The next result settles the question of the existence
(and non-existence) of solutions to the CYBE lying in
$\frak p_i\bigwedge \frak p_i$.
\proclaim{Theorem 5.5} Let $\frak p_i$ be the maximal parabolic
subalgebra of $\frak {sl}(n)$ obtained by deleting the
negative simple root $e_{i+1,i}$.
\roster
\item $\frak p_i$ is Frobenius if and only if $i$ and $n$ are
relatively prime.
\item $H^2(\frak p_i)=0$ for all $i$.
\endroster
\endproclaim
\demo{Proof} (1) This is a special case of the theorem in [E1]
which classifies Frobenius Lie algebras of the form
$R+N$, where $R$ is a reductive subalgebra and $N$ is a unipotent
radical which is either a simple $R$-module or abelian.
\newline
(2) This is an immediate consequence of
the Hochschild-Serre spectral sequence, see [F]. The result can
also be obtained directly by considering the decomposition
of $\frak p_i$ into its reductive and unipotent components. \qed
\enddemo
\proclaim{Corollary 5.6} A maximal parabolic subalgebra $\frak {sl}(n)$
is quasi-Frobenius if and only if it is Frobenius.
\endproclaim
\demo{Proof} Combine Theorem 5.5 with Definition 3.1. \qed
\enddemo
The preceeding two results imply that $\frak p_i$
is a carrier of a solution to the CYBE if and only if
$i$ and $n$ are relatively prime and, in this case, the
solution is uniquely determined up to equivalence.
These results however only establish the
existence of these solutions.
If $i$ and $n$ are relatively prime,
it seems natural to ask whether the solution to the CYBE with carrier
$\frak p_i$ is a boundary solution and, if so, is it possible
to explicitly construct it from an appropriate solution to the
MCYBE? We think that this is indeed the case. Recall that the
other instance where we needed $i$ and $n$ to be relatively
prime was in Theorem 2.6 which classified and constructed
all possible generalized
Cremmer-Gervais triples, $\Cal T_i$, for $\frak{sl}(n)$. Each of these
triples determines
a unique (up to equivalence) solution, $r_i$, to the MCYBE.
Now the deleted root from the subset
$\Pi_2$ of $\Pi$ associated to
$\Cal T_i$ is $i$.
Thus the maximal parabolic subalgebra $\frak p_i$
is, in a sense, canonically related
to both the unique solution to the CYBE (since it is Frobenius) and
the generalized Cremmer-Gervais triple $\Cal T_i$.
This prompts the following:
\definition{Conjecture 5.7} Suppose that $i$ and $n$ are relatively
prime and let $r_i$ be the unique solution to the MCYBE associated
to the generalized Cremmer-Gervais
triple $\Cal T_i$.
We conjecture that
the unique solution of the CYBE
with carrier $\frak p_i$ is a boundary solution and
lies in the closure of a suitable
orbit of $r_i$ under the action of $SL(n)$.
\enddefinition

We are able to verify this conjecture when $i=1$, (and hence also
for $i=n-1$) and also for the triple
$\Cal T_2$ associated to ${\frak {sl}}(5)$.
Our analysis of these cases suggests a procedure which we believe
works in general, see Conjecture 6.1.

We begin with $i=1$. This case corresponds to the
Cremmer-Gervais triple $\Cal T_{CG}$. The explicit form for the
associated solution, $r_{CG}$, to
the MCYBE is
$$\left(\sum_{i<j}e_{ij}\wedge e_{ji}\right)+
\left(\frac1n\sum_{i<j}(n+2(i-j))e_{ii}\wedge e_{jj}\right) +
2\left( \sum_{i<j}\sum_{m=1}^{j-i-1}e_{i,j-m}\wedge  
e_{j,i+m}\right).\tag5.8 $$
\newline\newline

\proclaim{Theorem 5.9} Suppose that
$r_{CG}\in \frak {sl}(n)\bigwedge \frak {sl}(n)$
is the Cremmer-Gervais solution to the MCYBE.
\roster
\item Let $x=\frac{1}{2}[(n-1)e_{12}+(n-2)e_{23}+\cdots
+1\cdot e_{n-1,n}]\in \frak {sl}(n)$ and set
$$d_p=\frac{n-p}{n}\left( e_{11}+e_{22}+\cdots e_{pp}\right)
- \frac{p}{n}\left( e_{p+1,p+1}+e_{p+2,p+2}+\cdots e_{nn}\right). $$
Then $[x,[x,r]]=0$ and so
$\exp (-t\, x)\cdot r_{CG}=r_{CG}+t\, [x,r]$. Thus 
$$[x,r]=b_{CG}=\left( \sum_{p=1}^{n-1}d_p\wedge e_{p,p+1}\right)
+\sum_{i<j} \sum_{m=1}^{j-i-1}e_{i,j-m+1}\wedge e_{j,i+m}$$
is a boundary solution to the CYBE.
\item The carrier for $b_{CG}$ is
$\frak p_1$, the maximal parabolic subalgebra of
$\frak {sl}(n)$ obtained by deleting the first negative root.
\item The inverse of $b_{CG}$ is the skew (non-degenerate)
bilinear form on $\frak p_1$ associated the the linear
functional $f=e_{12}^*+e_{23}^*+\cdots e_{n-1,n}^*:\frak  
p_1\rightarrow k$.
That is,
$b_{CG}^{-1}(p,q)=f([p,q])$ for all $p,q\in \frak p_1$.
\endroster
\endproclaim
\demo{Proof} Compute. \qed \enddemo

The theorem presents several
intriguing questions. 
First, notice that there is a striking similarity between the individual terms
of the last
summands of $r_{CG}$ and $b_{CG}$. For every
term $e_{i,j-m}\wedge e_{j,i+m}$ of the last summand of
$r_{CG}$, there is a similar term of
of the last summand of $b_{CG}$, namely $e_{i,j-m+1}\wedge e_{j,i+m}$.
The only difference is the $+1$ in the second subscript of the first
wedge factor. Is this related to the bijection $T$ from Cremmer-Gervais triple
$\Cal T_{CG}$? Recall that in this case for every
$i\in \Pi_1$ we had $T(i)=i+1$.

Theorem 5.9 also provides an interesting connection between
$r_{CG}$, $b_{CG}$, and the principal three-dimensional subalgebra
of $\frak {sl}(n)$. Recall that this subalgebra (which is isomorphic
to $\frak {sl}(2)$) has generators \newline
$E=(n-1)e_{12}+(n-2)e_{23}+\cdots +
1\cdot e_{n-1,n}$, \newline
$F=1\cdot e_{21}+2\cdot e_{23}+\cdots +
(n-1)e_{n,n-1}$,\quad and \newline
$H=(n-1)e_{11}+(n-3)e_{22}+\cdots + (3-n)e_{n-1,n-1}+(1-n)e_{nn}$. \newline
This copy of $\frak {sl}(2)$ naturally acts on
$\frak {sl}(n)\bigwedge \frak {sl}(n)$ via the adjoint action.
In terms of this action, the first part of Theorem 5.9 says that
$[E,r_{CG}]= -2b_{CG}$ and $[E,b_{CG}]=0$. We can relate this
to the $\frak {sl}(2)$-module structure
on $\frak {sl}(n)\bigwedge \frak {sl}(n)$ by
means of the automorphism
$\sigma : \frak {sl}(n)\rightarrow \frak {sl}(n)$ 
which sends $e_{ij}$ to $e_{n+1-i,n+1-j}$. It is easy
to see that $\sigma^2$ is the identity, 
$\sigma (E)=F,$ $\sigma (F)=E$, and
$\sigma (r_{CG})=-r_{CG}$. Therefore it follows that
$[F,[F,\sigma (b_{CG})]]=0$ and so $\sigma (b_{CG})=2[F,r_{CG}]$ is also
a boundary solution. Moreover, the three-dimensional submodule of
$\frak {sl}(n)\bigwedge \frak {sl}(n)$ generated by $r_{CG}$ under the action
of the principal three-dimensional subalgebra of ${\frak {sl}}(n)$ is
simple
since $\frak{sl}(n)$ decomposes into a direct sum of
simple submodules of dimensions $3, 5, 7, \dots,(2n-1)$ under this
action and hence
$\frak{sl}(n) \bigwedge \frak{sl}(n)$ contains no invariants (one-dimensional 
submodules). But the only possible decomposition of
the submodule generated by $r_{CG}$  would be into a simple three-dimensional
submodule and a direct sum of invariants; there being none of the latter, the
module must be simple. For this simple module, the highest and lowest weight
vectors, $b_{CG}$ and $\sigma(b_{CG})$ are boundary solutions to the CYBE
while the weight zero vector, $r_{CG}$, is a solution to the MCYBE.

Let us now consider the triple $\Cal T_2$ for $\frak {sl}(5)$
and its associated solution $r_2$ to the MCYBE. 
Explicitly, $r_2 =\gamma +\beta +\alpha $
where 
$$\split 
&\\
\gamma &= e_{12}\wedge e_{21}+e_{13}\wedge e_{31}+e_{14}\wedge e_{41}
+e_{15}\wedge e_{51} +e_{23}\wedge e_{32}\\
& \quad +e_{24}\wedge e_{42}+e_{25}\wedge e_{52}+e_{34}\wedge e_{43}
+e_{35}\wedge e_{53} +e_{45}\wedge e_{54},\\
&\\
\beta &=\frac {1}{5}\big( -e_{11}\wedge e_{22} -e_{22}\wedge e_{33} -
e_{33}\wedge e_{44} - e_{44}\wedge e_{55} + 3 e_{11}\wedge e_{33}\\
& \quad +3e_{22}\wedge e_{44} +3 e_{33}\wedge e_{55} -
3e_{11}\wedge e_{44} - 3e_{22}\wedge e_{55} +  e_{11}\wedge e_{55}\big),
\quad {\text {and}}\\
&\\
\alpha &=2( e_{23}\wedge e_{54} +e_{23}\wedge e_{21} +
e_{23}\wedge e_{43} + e_{45}\wedge e_{21}
+ e_{45}\wedge e_{43}+e_{12}\wedge e_{43}+e_{13}\wedge e_{53}).\\
&
\endsplit $$
To show that
the solution to the CYBE with carrier $\frak p_2$ lies in the boundary of
an orbit of $r_2$
we again make use of a three-dimensional subalgebra of $\frak {sl}(5)$, but
this one is not the principal one and the submodule of 
$\frak{sl}(5) \bigwedge \frak{sl}(5)$ generated by $r_2$ will not be simple.

Set $X = 2e_{13}+e_{24}+e_{35}, \quad Y = e_{31}+e_{42}+2e_{53}$ and $H
= [X,Y] = 2e_{11}+e_{22}-e_{44}-2e_{55}.$  Then we have $[H,X] =
2X,\quad [H,Y] = -2Y,$ and so the subalgebra of $\frak {sl}(5)$ 
spanned by $H, X$ and
$Y$ is isomorphic to $\frak{sl}(2)$. Let $\sigma: {\frak {sl}}(5)\rightarrow
\frak{sl}(5)$ once again denote the
automorphism sending $e_{i,j}$ to
$e_{n+1-i,n+1-j}$. Note that 
$\sigma (X) = Y, \quad \sigma( H) = H,$ and $\sigma (r_2) = 
-r_2.$  Just as in the proof of Theorem 5.9,
it is an easy computation to show that we have
$\exp (-tX)\cdot r_2 = r_2 + t[X,r_2]$ where
$$\split [X,r_2]& =
{\tsize {\frac{1}{5}}}
\{2(2e_{11}+2e_{22}-3e_{33}+2e_{44}-3e_{55})\wedge e_{13}+
(e_{11}+e_{22}+e_{33}-4e_{44}+e_{55})\wedge e_{24}\\
&\quad +
(e_{11}+e_{22}+e_{33}+e_{44}-4e_{55})\wedge e_{35}\}\\
&\quad+(e_{14}\wedge e_{43} +e_{12}\wedge e_{23} +
e_{25}\wedge e_{54} + e_{12}\wedge e_{45}+
e_{15}\wedge e_{53} + e_{34}\wedge e_{45})\endsplit $$
and so $[X,r_2]$ is a
boundary solution to the CYBE. 
However, its carrier is not
$\frak p_2,$ which has dimension 18, but a 16-dimensional
subalgebra. This is probably related to the fact that the present
[$X,r]$ is already in the boundary of the orbit of a ``smaller"
solution $r_2'$ to the MCYBE. 
 
Under the action of the current ${\frak {sl}}(2)$, we have that
$\frak{sl}(5)$ decomposes into a
direct sum of modules of dimensions $5, 4, 4, 3, 3, 2, 2,$ and
$1,$ so $\frak{sl}(n) \bigwedge \frak{sl}(n)$ has two invariants
(coming from the exterior products of the 4 dimensional and 2
dimensional submodules with themselves). 
Setting $\xi = e_{23}+e_{45}$ and $\eta = \sigma(\xi )=
e_{21}+e_{43},$ we can write $\alpha = \alpha_0 + \alpha_1$ where
$\alpha_0 = 2 \xi \wedge \eta$ and 
$\alpha_1 = 2(e_{23}\wedge e_{54}+e_{12}\wedge e_{43}+
e_{13}\wedge e_{53}).$  Then, with the present $X, Y, H,$ one can
check the following relations: $[X,\xi] = 0, \quad [X,\eta] = -\xi,$ and so
$[X,\alpha_0] = 0$. Similar relations for
$Y$ are easily obtained by applying $\sigma$. Thus
$\alpha_0$ is invariant for the action of the current copy of
$\frak {sl}(2)$ and so if
$r_2' = r_2 - \alpha_0$ then $[X,r_2] = [X,r_2']$ and $[X,r_2]$ is in the
boundary of the orbit of $r_2'$. It follows that 
the module
generated by $r_2'$ is 3-dimensional and simple, and the module
generated by $r_2$ is the direct sum of this module and the
invariant $\alpha_0.$ 
Now one can check that $r_2'$ is in fact
one of the Belavin--Drinfel'd solutions to the MCYBE, namely that
obtained by taking $\Pi_1 = {1,2}, \Pi_2 = {3,4},$ with $T ( 1 )=
3, T (2 )= 4.$ Note that this triple is obtained from $\Cal T_2$
by dropping the unique map, $T(4)=1$, which takes a root higher
in the usual order to a lower one.

Thus the direct analog of Theorem 5.9 does not produce the solution
to the CYBE with carrier $\frak p_2$. However, we can 
successfully modify the proceeding to indeed find this solution.
Let $H_1=(1/5)(-4e_{11}+6e_{22}-4e_{33}+6e_{44}-4e_{55})$.
Then $\exp (\eta)\cdot r_2= r_2+H_1\wedge \eta$ is a solution
to the MCYBE equivalent to $r_2$.
Now
$$\exp (tX+\eta)\cdot r_2 =r_2+H_1\wedge \eta +t\left\{
[X,r_2]-(3/2)H_1\wedge \xi +\eta \wedge \xi\right\}.$$
The coefficient of $t$ must therefore be a boundary solution and
one can check that its carrier is the maximal parabolic subalgebra
$$\frak p_2=\pmatrix *&*&*&*&*\\
*&*&*&*&*\\
0&0&*&*&*\\
0&0&*&*&*\\
0&0&*&*&*\endpmatrix$$

\subhead{6. Closing remarks}
\endsubhead
We close with two conjectures -- one about
how the foregoing analysis may extend to prove
Conjecture 5.7 in all cases and the other about constructing
quantum Yang-Baxter matrices from the solutions to the CYBE.

Let us consider the $i$th maximal parabolic subalgebra of $\frak{sl}(n)$
where $i$ is relatively prime to $n$ and hence defines a Belavin-
Drinfel'd solution, $r_i$, to the MCYBE in which the sole root omitted
from $\Pi_1$ is $n-i.$ As mentioned in Remark 2.7, the integers 
$1,\dots,n-1$ have an order defined by the corresponding root
mapping $T.$ That is, write them in the order $i, 2i, \dots, n-
i,$ always reducing modulo $n.$  So, for example, if $n=12$ and
$i=5$ then the order is $5,10,3,8,1,6,11,4,9,2,7.$  Now define a
{\it {string}}, $s,$ to be a maximal subsequence in which the integers
appear in their natural order. In this case, we have the
following strings: $\{5, 10\}, \{3,8\}, \{1,6,11\}, \{4,9\},
\{2,7\},$ which we will denote by $s_1,\dots,s_5,$ respectively,
and in the general case by $s_1,\dots,s_i$ since the number of
strings is always equal to $i.$ Some will have length $[n/i]$
(greatest integer contained in $n/i$) and some will have length
one greater.  Now for each string $s_1,\dots,s_{i-1},$ i.e., for
all except the last, define $e(s) = \sum e_{j,j+1}$ where the sum
runs over the integers $j$ in the string.  In the example we have
$e(s_1)=  e(5,10) = e_{5,6} + e_{10,11},$ but $e(2,7)$ is not
defined. For every string $s_2,\dots,s_i,$ i.e. all but the
first, define $e'(s) = \sum e_{j+1,j}$ where the sum again is
over the integers in the string.  Thus $e'(s_5) = e'(2,7) =
e_{3,2} + e_{8,7},$ while $e'(5,10)$ is not defined. Notice that
if, as before, $\sigma$ is the automorphism of $\frak{sl}(n)$ sending
$e_{jk}$ to $e_{n+1-j,n+1-k},$ then $\sigma$ sends the set of
$e(s)$ which are defined to the set of $e'(s)$ which are defined.
 
Let $L_+$ now be the Lie subalgebra of $\frak{sl}(n)$ generated by
$e(s_1),\dots,e(s_{i-1}),$ let $L_- = \sigma (L_+)$ be that
generated by $e'(s_2),\dots,e'(s_i),$ and $L$ be the Lie
subalgebra generated by both together. The intersection of $L$
with the set of diagonal matrices will be denoted $L_0.$  In the
case where $i = 1$ we have $L = 0,$ since  there is then but one
string and hence no $e(s)$ and no $e'(s)$ are defined.  The
algebra $L$ is never semisimple.  When $i=2$ it is of dimension 2
and abelian, but there is a non-trivial semisimple part whenever
$i>2.$  (In the case $n=5, i =2,$ the unique $e(s) = e(2,4)$ was
denoted by $\xi$ and the unique $e'(s) = e(1,3)$ by $\eta.$) The
radical $\Cal R$ of $L$ is the direct sum of its part lying in
$L_+,$ which we denote by $\Cal R_+,$ and of its part in $L_-,$
denoted $\Cal R_-.$  One can
readily compute that in the example $n=11, i=5$, the
semisimple part is a direct sum of two algebras, one isomorphic
to $\frak {sl}({2})$ and the other to $\frak {sl}({3}),$ 
its total dimension being 11,
while the radical has dimension 12.
Finally, we define a subalgebra isomorphic
to $\frak{sl}(2)$ as follows: set 
$X = \sum_{j=1}^{n-i}[(n-j)]/i] e_{j,j+i},\quad Y =
 \sigma X,$ and $H = [X,Y].$  
\definition{Conjecture 6.1}
\roster
\item $ [X,L_0] = [X,L_+] = 0$, while
$[X,L_-] = L_+$ and $[X, \Cal R_-] = \Cal R_+.$
\item $\Cal R$ is
abelian.  (One obtains the corresponding relations involving $Y$
by applying $\sigma.$) These assertions can be easily verified in
the example and probably are not too difficult in general.
The more difficult part is: 
\item Let
$r_i$ be the Belavin-Drinfel'd solution to the MCYBE associated to the
triple $\Cal T_i$ and suppose $z$ is a generic element of $\Cal R_-$ Then
$$\exp (tX+z)\cdot r_i =r_i' +t\omega$$
where $r_i'$ is a solution to the MCYBE which is equivalent to $r_i$
and $\omega $ is a boundary solution whose carrier is
precisely the maximal parabolic subalgebra $\frak p_i$.
(This is
essentially what we have shown for $n=5, p=2.$)
\endroster
\enddefinition
The Lie algebra $L$ arises naturally by considering the terms which 
disappear from the $\alpha$-part of $r_i$ when one drops from
the bijection $T$ of $\Cal T_i$ the mappings carrying a root higher
in the usual order to a lower one. But the significance of its
decomposition is not clear, and what is most mysterious is the 
representation of $\frak {sl}(2)$ associated with each $i$ which seems
to play an essential role.

We do not know in general which non-maximal parabolic subalgebras
of $\frak {sl}(n)$ are Frobenius or quasi-Frobenius, but if the conjecture
is true it may suggest ways to determine when a parabolic is the
carrier of a boundary solution to the CYBE (and hence at least
quasi-Frobenius).

Finally, it is natural to consider the problem of explicitly finding a quantum
Yang-Baxter matrix associated to a solution of the CYBE. Little is
known in general about this problem. For many of the boundary
solutions though the following result applies.
\proclaim{Theorem 6.2 [GGS1]} Suppose that
$b\in \frak{sl}(n)\bigwedge \frak{sl}(n)$ is a solution to the 
CYBE and that $b^3=0$ when viewed as an $n^2 \times n^2$ matrix under the
Kronecker product. Then $B=\exp (t\, b)=1 +t\,b +(t^2/2)b^2$ is a
solution to the quantum Yang-Baxter equation, that is,
$B_{12}B_{13}B_{23}=B_{23}B_{13}B_{12}$.
\endproclaim
The boundary solution of Example 5.2 has cube zero and it seems
likely that many others do as well. This is true of the
Cremmer-Gervais 
boundary solution $b_{CG}$ for $n=3$ (formula (0.2)). 
In this case we have that $1+t\,b_{CG}+(t^2/2)b^2_{CG} =$
$$\pmatrix
1&t/3&t^2/36&-t/3&t^2/18&0&t^2/36&0&0\\
\\
0&1&t/6&0&t/6&t^2/36&0&t^2/36&0\\
\\
0&0&1&0&0&t/6&0&t/2&-t^2/18\\
\\
0&0&0&1&-t/6&t^2/36&-t/6&t^2/36&0\\
\\
0&0&0&0&1&t/6&0&-t/6&t^2/18\\
\\
0&0&0&0&0&1&0&0&t/3\\
\\
0&0&0&0&0&-t/2&1&-t/6&-t^2/18\\
\\
0&0&0&0&0&0&0&1&-t/3\\
\\
0&0&0&0&0&0&0&0&1
\endpmatrix$$
is a solution to the quantum Yang-Baxter equation.
\proclaim{Conjecture 6.3} For all $n>3$, the Cremmer-Gervais boundary
solution $b_{CG}$ (Theorem 5.9.1) has cube zero.
\endproclaim

 \enddocument